\documentclass[12pt]{article}
\def\RA{\rightarrow}
\usepackage{theorem}
\usepackage{latexsym}
\usepackage{amssymb}
\usepackage{amsmath}
\usepackage{epic}
\usepackage{pifont}
\font\smallit=cmti10
\font\smalltt=cmtt10

\def\B{\hfill $\Box$}

\newtheorem{lemma}{Lemma}

\begin{document}
\thispagestyle{empty}

\vskip 30pt
\begin{center}
{\bf ON THE ASYMPTOTIC MINIMUM NUMBER OF MONOCHROMATIC
3-TERM ARITHMETIC PROGRESSIONS}
\vskip 20pt
{\bf Pablo A. Parrilo}\\
{\smallit  Department of Electrical Engineeering and Computer Science,\\
Massachusetts Institute of Technology,
Cambridge, MA 02139}\\
{\smalltt parrilo@mit.edu}
\vskip 10pt
{\bf Aaron Robertson}\\
{\smallit Department of Mathematics,
Colgate University,
Hamilton, NY 13346}\\
{\smalltt aaron@math.colgate.edu}
\vskip 10pt
{\bf Dan Saracino}\\
{\smallit Department of Mathematics,
Colgate University,
Hamilton, NY 13346}\\
{\smalltt dsaracino@mail.colgate.edu}

\end{center}
\vskip 30pt
\renewcommand{\arraystretch}{1.5}


\begin{abstract}
Let $V(n)$ be the minimum number of monochromatic
$3$-term arithmetic progressions in any $2$-coloring
of $\{1,2,\dots,n\}$.
We show that
$$\frac{1675}{32768} n^2(1+o(1)) \leq V(n) \leq
\frac{117}{2192}n^2(1+o(1)).$$ As a consequence, we find that $V(n)$
is strictly greater than the corresponding number for Schur triples
(which is $\frac{1}{22}n^2(1+o(1))$).  Additionally, we disprove the
conjecture that $V(n) = \frac{1}{16}n^2(1+o(1))$ as well as a more
general conjecture.
\end{abstract}

\pagestyle{myheadings}
\markright{\smallit On the Asymptotic Minimum Number of
Monochromatic Van der Waerden Triples}
\baselineskip=16pt

\section*{\normalsize 1. Introduction}

At the Erd\H{o}s Conference in Budapest in the summer of 1999,
Ron Graham proposed
the following \$100 problem:

\begin{quote}
{ Let $V(n)$ be the minimum number of monochromatic
$3$-term arithmetic progressions  in any
$2$-coloring of
$[1,n]=\{1,2,\dots,n\}$. Given $V(n) = \beta n^2(1+o(1))$,
determine $\beta$.}
\end{quote}

This problem seems to be much more abstruse than
the corresponding problem concerning Schur triples
(see [D], [RZ], [S]).  
It is conjectured, and commonly believed, that $\beta = \frac{1}{16}$,
in part because of the following ``folklore" conjecture.

\noindent
{\bf Conjecture.}   The minimum number of monochromatic
solutions, in any $r$-coloring of $[1,n]$, of
$
 \sum_{i=1}^m c_i x_i =0
$
with $\sum_{i=1}^m c_i=0$  is
equal to the value achieved by  randomly
coloring the integers in $[1,n]$.

In the case of $3$-term arithmetic progressions, the equation is
$
x+y=2z
$
and the value achieved by randomly $2$-coloring
the integers in $[1,n]$
is
$\frac{n^2}{16}(1+o(1))$ since there are $\frac{n^2}{4}(1+o(1))$
$3$-term arithmetic progressions in $[1,n]$, 
of which $\frac{1}{4}$ is the expected fraction of them
that are monochromatic under a random $2$-coloring.

The conjecture states that $V(n) = \frac{n^2}{16}(1+o(1))$.
We  show that this conjecture is  false
by proving that
$
V(n)<\frac{n^2}{16}(1+o(1)).
$
While we do not find $\beta$,
we are able to offer fairly good
upper and lower bounds.  We do believe
that the upper bound is extremely close,
if not equal, to $V(n)$.

\section*{\normalsize 2. Preliminaries for the Lower Bound}
\label{sec:prelim}

Let $\chi:[1,n] \RA \{0,1\}$ be an
arbitrary $2$-coloring.
Define, for $j=0,1$,  
$$S_j = \{x : \chi(x)=j, \; \; 1 \leq x \leq n\}.$$ 
Let $V(S_0,S_1)=V(n;S_0,S_1)$ be the number of monochromatic
$3$-term arithmetic progressions in $[1,n]$ under $\chi$.

Using an approach found in [S] and [D], we let
$$
f_j = \sum_{s \in S_j} e^{2 \pi i s x}, \,\,\,\, j=0,1,
$$
which gives us
$$
2V(S_0,S_1) = \int_0^1 \left(f^2_0(x) \overline{f_0(2x)} +
f^2_1(x) \overline{f_1(2x)}\right)\, dx.
$$

We rewrite the integrand as
$$
\begin{array}{rl} \displaystyle
\left(f_0(x)+f_1(x) \right)^2
\left(\overline{f_0(2x)}+\overline{f_1(2x)}\right)
&-
\left(f_0(x)\overline{f_1(2x)} + f_1(x)\overline{f_0(2x)}\right)
\left(f_0(x) + f_1(x)\right)\\  \displaystyle
&- \,\, f_0(x)f_1(x)
\left(\overline{f_0(2x)}+\overline{f_1(2x)}\right)
\end{array}
$$
and interpret the integral as
$$
\begin{array}{rl}
2 V(S_0,S_1) &= \left|\{(a,b,c) \in [1,n]^3: a+b=2c\}\right| \\ 
& \displaystyle - \left|\left\{(a,b) \in (S_0 \times S_1) \cup (S_1 \times
S_0): 2b-a
\in [1,n]\right\}
\right| \\ 
&  \displaystyle - \left| \left\{(a,b) \in S_0 \times S_1: a+b \mbox{ is
even}\right\}\right|.
\end{array}
$$

We will now bound the size of these sets, where our
equations are valid up to $o(n^2)$.

It is trivial to show that
$\left|\{(a,b,c) \in [1,n]^3: a+b=2c\}\right| = \frac{n^2}{2}(1+o(1))$.
It is also easy to show that
$\left| \left\{(a,b) \in S_0 \times S_1: a+b \mbox{ is
even}\right\}\right| \leq \frac{n^2}{8}(1+o(1))$ as follows.
Denote this set by $T$ and
let $r_{o}$ and $b_o$ be the number of odd numbers in $[1,n]$
of color red (in $S_0$, say) and blue (in $S_1$), respectively,
 and let $r_{e}$ and $b_e$ the number of even numbers in $[1,n]$
of color red and blue, respectively.  
Then 
$$
\begin{array}{rl}
|T|&=(r_ob_o + r_e b_e )\\
&= \frac12\left((r_o+b_o)^2+(r_e+b_e)^2 - (r_o^2+r_e^2+b_o^2+b_e^2)\right)\\
&= \frac12\left(\left(\frac{n}{2} \right)^2+ \left(\frac{n}{2} \right)^2
- (r_o^2+b_o^2+r_e^2+b_e^2)\right)\\
&=\frac12\left(\frac{n^2}{2}-
(r_o^2+b_o^2+r_e^2+b_e^2)\right)\\
&=\frac{n^2}{4} -
\frac12\left(r_o^2+\left(\frac{n}{2}-r_o\right)^2
+r_e^2+\left(\frac{n}{2}-r_e\right)^2\right) \\
&= \frac{n}{2}(r_o+r_e) - (r_o^2+r_e^2).\\
\end{array}
$$
This function attains its maximum of
$\frac{n^2}{8}(1+o(1))$ when $r_o=r_e= \frac{n}{4}.$

Next, we define
$$N^+ = \left\{(a,b) \in (S_0 \times S_1) \cup (S_1 \times
S_0): 2b-a
\in [1,n]\right\}.$$
Our goal is to find an upper bound for $|N^+|$ and use the
following lemma, which follows immediately from the paragraphs above.

\begin{lemma}
If $|N^+|\leq cn^2(1+o(1))$, then
$$
V(S_0,S_1) \geq \frac12\left(\frac{3}{8}-c\right) n^2(1+o(1)).
$$
\end{lemma}

\section*{\normalsize 3. Lower Bound Calculations}
\label{sec:lower}

Our approach will be to consider points in the square $[1,n]^2$.  From
the definition of $N^+$, we restrict our attention to those points
$(x,y)$ with $0<2y-x \leq n$.  We also remark that since we are
looking for the coefficient of the $n^2$ term in $V(n)$, we will
disregard points that contribute $o(n^2)$ to $V(n)$.


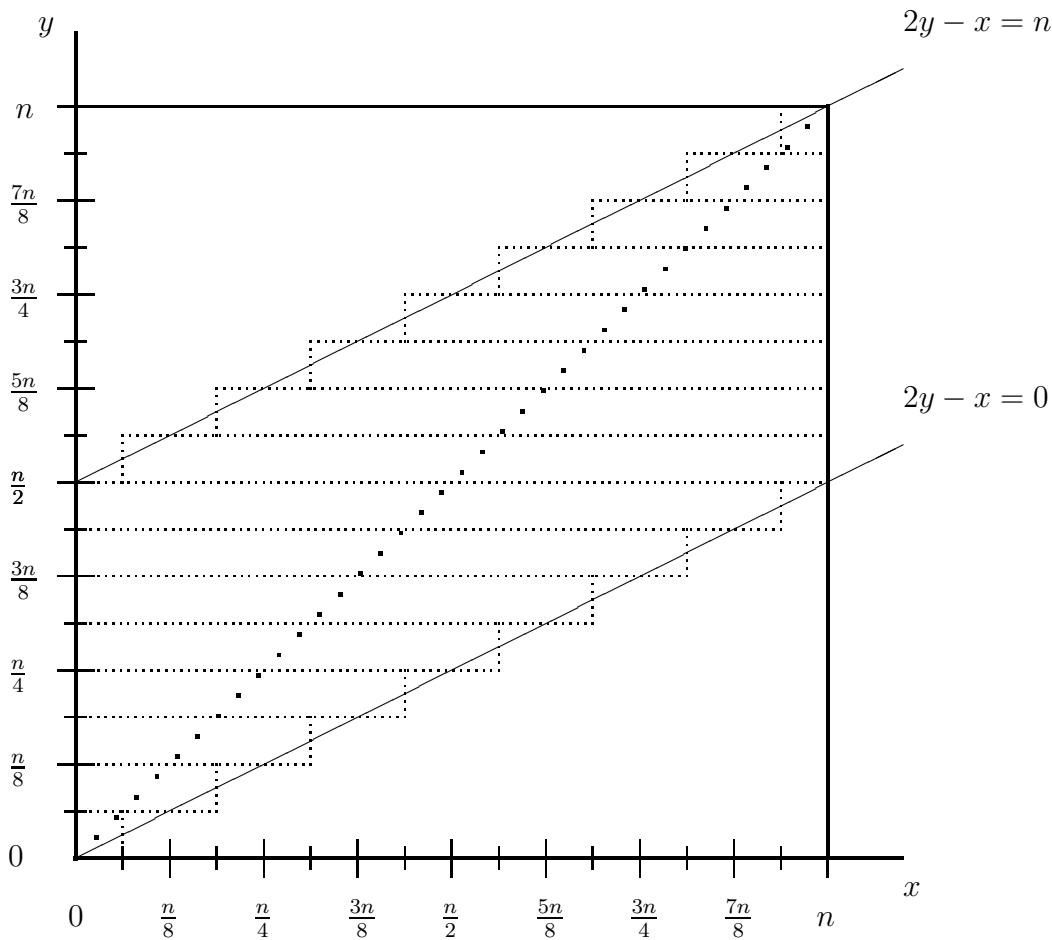
\begin{figure}
\setlength{\unitlength}{1mm}
\begin{picture}(0,130)(-30,-10)
\linethickness{.3mm}
\put(0,0){\line(0,1){110}}
\put(0,0){\line(1,0){110}}
\put(100,0){\line(0,1){100}}
\put(0,100){\line(1,0){100}}

\linethickness{.2mm}

\put(-2.5,12.5){\line(1,0){5}}
\put(-2.5,25){\line(1,0){5}}
\put(-2.5,37.5){\line(1,0){5}}
\put(-2.5,50){\line(1,0){5}}
\put(-2.5,62.5){\line(1,0){5}}
\put(-2.5,75){\line(1,0){5}}
\put(-2.5,87.5){\line(1,0){5}}
\put(-2.5,100){\line(1,0){5}}

\put(6.25,-1.5){\line(0,1){3}}
\put(18.75,-1.5){\line(0,1){3}}
\put(31.25,-1.5){\line(0,1){3}}
\put(43.75,-1.5){\line(0,1){3}}
\put(56.25,-1.5){\line(0,1){3}}
\put(68.75,-1.5){\line(0,1){3}}
\put(81.25,-1.5){\line(0,1){3}}
\put(93.75,-1.5){\line(0,1){3}}

\put(-1.5,6.25){\line(1,0){3}}
\put(-1.5,18.75){\line(1,0){3}}
\put(-1.5,31.25){\line(1,0){3}}
\put(-1.5,43.75){\line(1,0){3}}
\put(-1.5,56.25){\line(1,0){3}}
\put(-1.5,68.75){\line(1,0){3}}
\put(-1.5,81.25){\line(1,0){3}}
\put(-1.5,93.75){\line(1,0){3}}

\put(12.5,-2.5){\line(0,1){5}}
\put(25,-2.5){\line(0,1){5}}
\put(37.5,-2.5){\line(0,1){5}}
\put(50,-2.5){\line(0,1){5}}
\put(62.5,-2.5){\line(0,1){5}}
\put(75,-2.5){\line(0,1){5}}
\put(87.5,-2.5){\line(0,1){5}}
\put(100,-2.5){\line(0,1){5}}

\linethickness{.1mm}
\put(-9,-1){0}
\put(-1,-9){0}
\put(-5,110){$y$}
\put(110,-5){$x$}

\put(-9,48.5){$\frac{n}{2}$}
\put(-8,98.5){$n$}

\put(11,-9){$\frac{n}{8}$}
\put(23.5,-9){$\frac{n}{4}$}
\put(36,-9){$\frac{3n}{8}$}
\put(48.5,-9){$\frac{n}{2}$}
\put(61,-9){$\frac{5n}{8}$}
\put(73.5,-9){$\frac{3n}{4}$}
\put(86,-9){$\frac{7n}{8}$}

\put(-9,11){$\frac{n}{8}$}
\put(-9,23.5){$\frac{n}{4}$}
\put(-9,36){$\frac{3n}{8}$}
\put(-9,48.5){$\frac{n}{2}$}
\put(-9,61){$\frac{5n}{8}$}
\put(-9,73.5){$\frac{3n}{4}$}
\put(-9,86){$\frac{7n}{8}$}

\put(98.5,-9){$n$}

\linethickness{.2mm}

\dottedline{1}(0,50)(100,50)
\dottedline{1}(0,43.75)(93.75,43.75)
\dottedline{1}(0,37.5)(81.25,37.5)
\dottedline{1}(0,31.25)(68.75,31.25)
\dottedline{1}(0,25)(56.25,25)
\dottedline{1}(0,18.75)(43.75,18.75)
\dottedline{1}(0,12.5)(31.25,12.5)
\dottedline{1}(0,6.25)(18.75,6.25)

\dottedline{1}(6.25,56.25)(100,56.25)
\dottedline{1}(18.75,62.5)(100,62.5)
\dottedline{1}(31.25,68.75)(100,68.75)
\dottedline{1}(43.75,75)(100,75)
\dottedline{1}(56.25,81.25)(100,81.25)
\dottedline{1}(68.75,87.5)(100,87.5)
\dottedline{1}(81.25,93.75)(100,93.75)

\dottedline{1}(93.75,50)(93.75,43.75)
\dottedline{1}(81.25,43.75)(81.25,37.5)
\dottedline{1}(68.75,37.5)(68.75,31.25)
\dottedline{1}(56.25,31.25)(56.25,25)
\dottedline{1}(43.75,25)(43.75,18.75)
\dottedline{1}(31.25,18.75)(31.25,12.5)
\dottedline{1}(18.75,12.5)(18.75,6.25)
\dottedline{1}(6.25,6.25)(6.25,0)

\dottedline{1}(93.75,100)(93.75,93.75)
\dottedline{1}(81.25,93.75)(81.25,87.5)
\dottedline{1}(68.75,87.5)(68.75,81.25)
\dottedline{1}(56.25,81.25)(56.25,75)
\dottedline{1}(43.75,75)(43.75,68.75)
\dottedline{1}(31.25,68.75)(31.25,62.5)
\dottedline{1}(18.75,62.5)(18.75,56.25)
\dottedline{1}(6.25,56.25)(6.25,50)


\put(110,60){$2y-x=0$}
\put(110,110){$2y-x=n$}

\linethickness{.5mm}

\put(0,0){\line(2,1){110}}
\put(0,50){\line(2,1){110}}
\put(0,0){\dottedline{4}(0,0)(100,100)}



\end{picture}
\caption{Partition of the square into rectangles, for $L= 16$.}
\label{fig:partition}
\end{figure}

Consider the diagram in Figure~\ref{fig:partition}.  We are trying to
find the maximum number of dichromatic pairs $(a,b)$ that can reside
inside the parallelogram bounded by the lines $x=0$, $x=n$, $2y-x=0$,
and $2y-x=n$.  To this end, we cover the parallelogram by $L$
horizontal strips of height $\frac{n}{L}$ and right triangles with
dimensions $\frac{n}{2L} \times \frac{n}{L}$ (in
Figure~\ref{fig:partition}, we have $L=16$).  As such, we cover more
that the parallelogram (we have right triangles outside of the
parallelogram).  Hence, by maximizing the number of dichromatic pairs
inside the strips and the right triangles, we have an upper bound on
the maximum number of dichromatic pairs that can reside inside the
parallelogram.

Let $((i-1)\frac{n}{L},i\frac{n}{L}]$ contain $r_i$ red elements, $i =
1,2,\dots,L$.  Choosing $L$ to be even, we can easily write down a
formula for the number of dichromatic pairs that reside in the
horizontal strips:
\begin{equation}
\sum_{i=1}^{L/2} \sum_{j=1}^{2i-1} \left(r_i
\left(\frac{n}{L}-r_j\right)+\left(\frac{n}{L}-r_i\right) r_j\right) +
\sum_{i=L/2+1}^{L}\sum_{j=2i-L}^{L} \left(r_i \left(\frac{n}{L}-r_j\right)
+\left(\frac{n}{L}-r_i\right) r_j \right).
\end{equation}

What remains are the maximum possible number of dichromatic points in
the $L$ remaining triangles.  For these we use the trivial bound of
their areas, $L \times \frac{1}{2} \frac{n}{L} \frac{n}{2L}=
\frac{n^2}{4L}$. Combining this with (1), we have an upper bound on
$|N^+|$:

\begin{multline}
|N^+| \leq \frac{n^2}{4 L} + \displaystyle \sum_{i=1}^{L/2} \sum_{j=1}^{2i-1}
\left(r_i
\left(\frac{n}{L}-r_j\right)+\left(\frac{n}{L}-r_i\right) r_j\right) \\
\displaystyle+\sum_{i=L/2+1}^{L}\sum_{j=2i-L}^{L} \left(r_i \left(\frac{n}{L}-r_j\right)
+\left(\frac{n}{L}-r_i\right) r_j \right).
\label{eq:nplus}
\end{multline}

We present next two different techniques to effectively bound the
right-hand side of~(\ref{eq:nplus}). The first one relies on an
explicit enumeration of all the critical points (for $L=16$), while
the second approach uses a procedure based on semidefinite
programming.

\subsection*{\normalsize 3.1 Enumeration Bounds for $L=16$}
\label{sec:enum}

In this approach, all critical points in
$\left(0,\frac{n}{16}\right)^{16}$ are compared against all maximum
values at the $3^{16}-1$ boundary problems. The maximization problem has
been programmed into Maple as a small program called {\tt PABLO} and
the code is available from the second author's website\footnote{{\tt
http://math.colgate.edu/$\sim$aaron/programs.html}}.

After running for
approximately 136 hours on a 2.7GHz G5 Macintosh server, we find that
$$
|N^+| \leq \frac{579}{2048}n^2(1+o(1)).
$$
One  coloring that achieves this bound is
$$
(r_1,r_2,\dots,r_{16}) = \left(
\frac{7n}{128},\frac{7n}{128},0,\frac{7n}{128},\frac{n}{16},
0,0,0,\frac{n}{16},\frac{n}{16},\frac{n}{16},\frac{n}{128},
0,\frac{n}{16},\frac{n}{128},\frac{n}{128}
\right).$$
\normalsize

Applying Lemma 1, the above result gives us the following theorem.

\noindent
{\bf Theorem 2}  $V(n) \geq \frac{189}{4096}n^2(1+o(1))$.

\subsection*{\normalsize 3.2 Semidefinite Bounds}
\label{sec:semidef}

A different, more powerful way of bounding $|N^+|$ is based on
semidefinite relaxations. For this, consider first the change of
variables $r_i := \frac{1+x_i}{2} \frac{n}{L}$, so $r_i \in [0,
\frac{n}{L}]$ if and only if $x_i \in [-1,1]$. Then, equation
(\ref{eq:nplus}) can be written as
\[
\begin{split}
|N^+| & \leq 
\frac{n^2}{4 L} +
\displaystyle\sum_{i=1}^{L/2} \sum_{j=1}^{2i-1}
\frac{n^2}{2 L^2}(1-x_i x_j) 
\displaystyle+\sum_{i=L/2+1}^{L}\sum_{j=2i-L}^{L} 
\frac{n^2}{2 L^2}(1-x_i x_j) \\
& \leq 
\frac{n^2}{4 L} + \frac{n^2}{4} - \frac{n^2}{4 L^2} \, 
q(\mathbf{x}),
\end{split}
\]
where 
\[
q(\mathbf{x}) := 
\displaystyle\sum_{i=1}^{L/2} \sum_{j=1}^{2i-1}
2x_i x_j 
\displaystyle+\sum_{i=L/2+1}^{L}\sum_{j=2i-L}^{L} 
2x_i x_j.
\]
Our objective is to bound $|N^+|$ from above. For this, it is clearly
enough to obtain a lower bound of the quadratic form $q(\mathbf{x})$
over $[-1,1]^n$. This quadratic form can be represented as
$q(\mathbf{x}) = \mathbf{x}^T A \mathbf{x}$, where $A$ is an $L \times L$
symmetric integer matrix, with entries $A_{ij} = B_{ij}+B_{ji}$ and 
\[
B_{ij} = 
\begin{cases}
1 & \text{if } j +1 \leq 2 i \leq j + L \\
0 & \text{otherwise.} 
\end{cases}
\]
A useful bound for quadratic forms on the unit hypercube, used
extensively in the combinatorial optimization literature, can be
obtained as follows.

\noindent
{\bf Lemma 3}
Let $A$ be an $n \times n$ matrix and
let $D = \textrm{diag}(d_1, \ldots, d_n)$ be a diagonal matrix, such
that $A + D$ is positive semidefinite. Then, for all
$\mathbf{x} \in [-1,1]^n$,
$\mathbf{x}^T A \mathbf{x}$ is bounded below by $-\sum_{i=1}^n d_i$.

\noindent
{\it Proof.} Consider any vector $\mathbf{x} \in [-1,1]^n$. Since
$A+D$ is
  positive semidefinite it follows that
\[
0 \leq \mathbf{x}^T(A+D)\mathbf{x} = \mathbf{x}^T A \mathbf{x} +
\sum_{i=1}^n d_i x_i^2.
\]
Since $x_i^2 \leq 1$, we have $\mathbf{x}^T A \mathbf{x} \geq
-\sum_{i=1}^n d_i x_i^2 \geq -\sum_{i=1}^n d_i$. \B

For any finite value of $L$, a suitable set of $d_i$ can be found by
semidefinite programming. For the case $L=128$ we have found a
particular solution (given in the Appendix) using the SDP solver
SeDuMi, followed by a straightforward rounding procedure (to obtain
rational solutions). For such a solution, it can be easily verified on a
computer that the $128 \times 128$ rational matrix $A+D$ is indeed
positive definite.  Since we have $\sum_{i=1}^L d_i = 1364$, this
gives an upper bound for
$|N^+|$ with $c = \frac{1}{4 L}+\frac{1}{4} + \frac{1364}{4
L^2}=
\frac{4469}{16384}$, resulting in the lower bound 
(via Lemma 1) given in the next theorem.

\noindent
{\bf Theorem 4}
$V(n) \geq 
\frac{1675}{32768}n^2(1+o(1)).
$

\section*{\normalsize 4. The Upper Bound}

\noindent
{\bf Theorem 5}  $V(n) \leq \frac{117}{2192}n^2(1+o(1))$

\noindent
{\it Proof.}  Let $i^m = \underbrace{i i \dots i}_m$, i.e., a string of
$i$'s of length
$m$. Consider the coloring, using the colors $0$ and $1$,
\large
$$
0^{\frac{28}{548}n}\,\,1^{\frac{6}{548}n}\,\,0^{\frac{28}{548}n}
\,\,1^{\frac{37}{548}n}\,\,0^{\frac{59}{548}n}\,\,1^{\frac{116}{548}n}
\,\,0^{\frac{116}{548}n}\,\,1^{\frac{59}{548}n}
\,\,0^{\frac{37}{548}n}\,\,1^{\frac{28}{548}n}\,\,0^{\frac{6}{548}n}\,\,
1^{\frac{28}{548}n}.
$$
\normalsize
It is tedious -- but routine -- to show that under this coloring there are
$
\frac{117}{2192}n^2(1+o(1))
$ monochromatic $3$-term arithmetic progressions, thereby proving the theorem.
\B

The above coloring was found using a combination of computational and
analytic methods. We briefly describe these next.



As we have seen in the previous sections, the problem can be
essentially reduced to the minimization of the quadratic form
$q(\mathbf{x})$ over the unit hypercube. To understand the behavior of
the solution, we solved instances of this problem for large values of
$n$ ($n \approx 2000$). For this, a ``good'' initial candidate
coloring was found using the solution of the semidefinite relaxation,
followed by a randomization procedure known as Goemans-Williamson
rounding [GW]. The near-optimal solutions found all shared some nice
structural features, essentially being constant over large ranges of
$n$, with a small number of breakpoints (equal to 12 for most
solutions).

We then used a continuous approximation to the minimization of
$q(\mathbf{x}) = \mathbf{x}^T A \mathbf{x}$, given by
\[
\min_\phi \int_{-1}^1 \int_{-1}^1 k(x,y) \phi(x) \phi(y) dx dy, 
\]
where the function $\phi$ must satisfy $|\phi(\cdot)| = 1$ and the
kernel $k(x,y)$ is piecewise constant.  Based on the numerical
solutions for large $n$, we chose an ansatz where the function $\phi$
is symmetric ($\phi(x) = \phi(-x)$) and piecewise constant on 12
different intervals. 

Because $k(x,y)$ is piecewise constant, the objective function is a
\emph{piecewise quadratic} function of 5 variables, namely the
breakpoints (5 variables rather than 12 since we are
assuming symmetry). It turns out that, on the
partition associated with the solution obtained by numerical computation,
this function is strictly convex and its minimum lies inside the
partition. Solving for the (local) minimum of this quadratic function, we
obtained the breakpoints corresponding to the solution in Theorem~5. The
solution presented is thus ``locally optimal'' in the sense that no small
perturbation of the breakpoints will achieve a better value. Of
course, in principle the possibility remains that there exist
solutions of different structure that achieve even smaller values, so
the argument given is not enough to prove global optimality.

The Maple code that computes this quadratic function and performs the
minimization is also available at the location cited earlier.

\section*{\normalsize 5. Remarks for Further Investigation}

Clearly, the parallelogram described at the beginning of
Section 3 could be further refined by using larger
values of $L$.

For the enumeration technique in Section 3.1 
this would
provide sharper bounds, which converge to the optimal constant
$\beta$.  However, since the number of points to be checked grows
exponentially with $L$, there would be an enormous increase in the
computational cost (for example, adding two more variables would
increase the computing time to approximately 51 days). A possible
improvement here could be obtained by finding an upper bound on the
triangles for which we have used the trivial bound of their area,
although this would not help with the exponential behavior.

For the semidefinite bounds in Section 3.2, it is
relatively straightforward (and computationally feasible) to provide
slightly better lower bounds by increasing the value of $L$. However,
even if we let $L \rightarrow +\infty$, the obtained bounds will
likely \emph{not} converge to the optimal value of $\beta$, as there
seems to be an ``irreducible'' gap between the original problem and
its corresponding semidefinite relaxation. While this issue is
relatively well-understood for finite problems, it would be of
interest to fully understand the situation in this infinite limit.

Given our belief that the bound presented in Theorem~5 is sharp,
perhaps the most promising approach would be to attempt to directly
prove the (asymptotic) global optimality of the corresponding
solution.

\section*{\normalsize Acknowledgement}  We thank Ron Graham for
some background information, but more importantly for providing the
combinatorics community with interesting and challenging problems. The
first author
thanks Bernd Sturmfels for forwarding the email (from Doron
Zeilberger) that spurred this collaboration.

\section*{\normalsize References}
\footnotesize
[D] B. Datskovsky, On the number of monochromatic
Schur triples, 
{\it Advances in Applied Math.} {\bf 31} (2003), 193-198.
\vskip -5pt
\noindent
[GW] M. X. Goemans and D. P. Williamson,
{Improved Approximation Algorithms for Maximum Cut and Satisfiability Problems Using Semidefinite Programming},
{\it Journal of the ACM} {\bf 42(6)} (1995), 115-1145.
\vskip -5pt
\noindent
[GRR] R. Graham, V. R\"{o}dl, and A. Ruci\'nski,
{On Schur Properties of Random Subsets of Integers},
{\it J. Number Theory} {\bf 61} (1996), 388-408.
\vskip -5pt
\noindent
[RZ] A. Robertson, D. Zeilberger,
{A $2$-Coloring of $[1,N]$ Can Have $N^2/22+O(N)$
Monochromatic Schur Triples, But Not Less!},
{\it Elect. J. of Comb.} {\bf 5(1)} (1998), R19.
\vskip -5pt
\noindent 
[S] T. Schoen, { On the Number of Monochromatic Schur Triples}, 
{\it Europ. J. of Comb.} {\bf 20(8)} (1999), 855-866.

\section*{\normalsize Appendix}

A particular solution for the $d_i$ in Lemma 3 is
given by the numbers below.
\[
\begin{array}{cccccccccccccccccccccccccccccccccccccccccccccc} 
d = \frac{1}{4} \; [  
& 27&22&14&14&13&11&5&2&9&14&20&24&26&29&28&26\\
& 26&26&26&25&24&23&23&21& 22&21&27&30&37&41&48&50\\
& 54&53&53&53&53&55&59&65&70&76&79&83&84&86&84&81\\
& 74&69&61&53&49&50&56&61&66&65&61&51&46&46&41&37\\
& 37&41&46&46&51&61&65&66& 61&56&50&49&53&61&69&74\\
& 81&84&86&84&83&79&76&70&65&59&55&53&53&53&53&54\\
& 50&48&41&37&30&27&21&22&21&23&23&24&25&26&26&26\\
& 26&28&29&26&24&20&14&9& 2&5&11&13&14&14&22&27&].
\end{array} 
\]

\end{document}